# Two-Phase Framework for Optimal Multi-Target Lambert Rendezvous

Jun Bang[1] and Jaemyung Ahn[2]

*Korea Advanced Institute of Science and Technology (KAIST), 291 Daehak-Ro, Daejeon 34141, Republic of Korea*

## Abstract

This paper proposes a two-phase framework to solve an optimal multi-target Lambert rendezvous problem. The first phase solves a series of single-target rendezvous problems for all departure-arrival object pairs to generate the elementary solutions, which provides candidate rendezvous trajectories (elementary solutions). The second phase formulates a variant of traveling salesman problem (TSP) using the elementary solutions prepared in the first phase and determines the best rendezvous sequence and trajectories of the multi-target rendezvous problem. The validity of the proposed optimization framework is demonstrated through an asteroid exploration case study.

## Nomenclature

| | |
|---|---|
| $i$ / $j$ | indices representing departure/arrival objects for a rendezvous |
| $p$ | index representing an arc |

[1] Graduate Research Assistant, Department of Aerospace Engineering, 291 Daehak-Ro.

[2] Associate Professor of Aerospace Engineering, 291 Daehak-Ro; jaemyung.ahn@kaist.ac.kr (Corresponding Author).

| | |
|---|---|
| $A$ | set of arcs in a graph |
| $c$ | cost associated with an arc in a graph |
| $G$ | graph representing a multi-target problem |
| $J$ | objective function |
| $l$ | line of intersection for the single target rendezvous |
| $m_i$ | number of revolutions of $i^{th}$ rendezvous, $\boldsymbol{\mu} = [m_1, \ldots, m_N]$ |
| $N$ | number of targets |
| $n$ | mean motion of an orbit |
| $\mathbf{P}_{\mathrm{M}}$ | optimal multi-target Lambert rendezvous problem |
| $\mathbf{P}_{\mathrm{ME}}$ | optimal multi-target Lambert rendezvous problem with elementary solutions |
| $\mathbf{P}_{\mathrm{S}}$ | optimal single-target Lambert rendezvous problem |
| $Q$ | set of targets, $Q = \{1, \ldots, N\}$ |
| $q_i$ | index of $i^{th}$ rendezvous target, $q_i \in Q$ |
| $\mathbf{q}$ | sequence of multi-target rendezvous, $\mathbf{q} = [q_1, \ldots, q_N]$ |
| $\mathbf{r}$ / $\mathbf{v}$ | position/velocity vectors |
| $t_{1,i}$ / $t_{2,i}$ | time at departure/arrival of $i^{th}$ rendezvous, $\boldsymbol{\tau}_1 = [t_{1,1}, \ldots, t_{1,N}]$ |
| $t_{ser}$ | required service time on a target |
| $\Delta t_{tr,i}$ | transfer time of $i^{th}$ rendezvous, $\Delta\boldsymbol{\tau}_{tr} = [\Delta t_{tr,1}, \ldots, \Delta t_{tr,N}]$ |
| $V$ | set of nodes in a graph |
| $\mathbf{x}^E_{(i,j)}$ | elementary solutions for an object pair ($i$, $j$) |
| $y$ | decision variable for $\mathbf{P}_{\mathrm{ME}}$ |
| $\theta$ | angle between position vector and line of intersection |

# I. Introduction

A multi-target rendezvous problem determines the visiting sequence and associated trajectories of a spacecraft to rendezvous with multiple objects. The problem is applicable to a number of space mission categories such as active debris removal (ADR), on-orbit servicing, and interplanetary exploration, and thus has been attracting considerable attention recently. Unlike a single-target rendezvous problem, the optimal multi-target rendezvous problem determines two different variable types (continuous and discrete) simultaneously. The problem is categorized as a mixed integer nonlinear programming (MINLP), which is known as one of the most difficult problem classes.

This paper proposes a two-phase framework to solve the multi-target Lambert rendezvous problem. The first phase of the framework generates the “elementary solutions,” which are the candidate components of the final solution, by solving a series of single-target rendezvous problems for all departure-arrival object pairs. The second phase combines the elementary solutions prepared in the first phase to obtain the final solution of the original multi-target rendezvous problem – best rendezvous sequence and trajectories. In the framework, the elementary solutions (representing candidate rendezvous trajectories) are used to transform the multi-target rendezvous problem into a variant of traveling salesman problem (TSP) that has multiple arcs between each pair of nodes and time window constraints associated with each arc.

Key contribution of this work is threefold. First, we introduced a framework that decomposes the original multi-target rendezvous problem categorized as the MINLP into a series of nonlinear programming (NLP) and an integer linear programming (ILP). The proposed framework effectively reduces the complexity of trajectory optimization process required to obtain the final solution, and provides better solutions than known approaches for the problem (e.g. variants of genetic algorithm) with reasonable increase in computational resource consumption. Secondly, the resulting ILP is a routing problem class that we named as “TSP with multiple arcs and arc time windows.” This problem is distinguishable from existing routing problems that consider only multiple arcs between nodes (Garaix et al. 2010; Ticha et al. 2017) or only arc time windows (Cetikaya et al. 2013) in that it deals with both of them. Finally, realistic case studies that can demonstrate the effectiveness of the proposed framework have been conducted. The multiple asteroid rendezvous mission is selected as the subject of the case study. The solutions for the case problems are obtained using the proposed framework and compared with known results found with other algorithm.

The rest of this paper is organized as follows. Section II provides reviews on past studies about the optimal rendezvous problems. Section III presents the mathematical formulation of optimal multi-target Lambert rendezvous mission. The detailed explanation on the steps of the framework to solve the problem is provided in Section IV. Case studies that can demonstrate the effectiveness of the proposed framework is presented in Section V. Finally, Section VI discusses the conclusions and potential future work of this study.

## II. Literature Review

An optimal single-target rendezvous, which underlies the multi-target problem, has been addressed in a number of published studies. One of widely used approaches for this problem is the *Lambert rendezvous*, which uses the solution of the Lambert's boundary value problem to obtain the rendezvous trajectory (Battin 1999). Past studies demonstrated that allowing multiple-revolution solutions of the Lambert's problem could reduce the fuel consumption required for the rendezvous. When we allow $N_{max}$ for a given transfer time there exist ($2N_{max}$+1) trajectories, which enlarges the search space for the design variables characterizing the rendezvous (Prussing 2000). Shen and Tsiotras (2003) proposed an algorithm that can determine an optimal solution of the fixed-time, two-impulse rendezvous between coplanar circular orbits among ($2N_{max}$+1) solutions quickly and efficiently by solving the multiple-revolution Lambert's problem. Zhang et al. (2011) proposed a procedure to solve the optimal two-impulse rendezvous problem for non-coplanar elliptic orbits. They considered coasting as an additional mission element and used a genetic algorithm (GA) to determine the initial and final coasting periods that minimize the propellant consumption. Chen et al. (2013) proposed a time-open constrained Lambert rendezvous problem by introducing parking time and transfer time. An interval branch-and-bound algorithm was adopted in combination with a gradient-based algorithm to find the global solution of the problem while effectively dealing with its strong nonlinearity and nonconvexity.

Optimal multi-target rendezvous problems are composed of two main tasks: determining rendezvous sequence and rendezvous trajectory optimization. Most of previous studies on the multi-target rendezvous focused on exploring the visiting sequences while optimization of the rendezvous trajectories was not seriously addressed – assumption-driven simplistic strategies were adopted or trajectory related discussions were missing. Determination of a visiting sequence for the multi-target rendezvous is similar to, but much more difficult to handle than the TSP primarily because the nodes (targets) are moving and the cost spent on an arc (rendezvous maneuver) is not fixed. Alfriend et al. (2005) formulated a geosynchronous satellites servicing mission as a relaxed TSP considering the fuel consumption only for orbital plane changes – cost associated with in-plane maneuver is ignored. Izzo et al. (2015) also estimated the cost as a relative orbital inclination and introduced two types of TSP variants for an active debris removal (ADR) mission – considering both of static and dynamic cases. Note that the optimization of the rendezvous trajectories was not the primary focus of either of the studies. Barbee et al. (2012) proposed a series method that can find a good – while not necessarily optimal – visiting order for the ADR mission with relatively low computational load. Some researchers, including Cerf (2013; 2015), Bérend and Olive (2016), introduced a drift orbit, which is an intermediate orbit that can accelerate the right ascension of the ascending node (RAAN) drift of a spacecraft. A branch and bound (B&B) algorithm was then employed to solve the optimal sequencing problem for an ADR mission using the cost of a specific three-step transfer strategy based on the Hohmann transfer.

On the other hand, there are several studies that address the Lambert rendezvous based trajectory optimization and the optimal rendezvous sequence determination simultaneously. Zhang et al. (2014; 2015) formulated the problems arising in mission designs for multi-spacecraft refueling and asteroid exploration as the multi-target rendezvous, and optimized the visiting sequence and the rendezvous trajectories simultaneously by introducing a procedure based on the hybrid encoding genetic algorithm (HEGA). Ross and D’Souza (2005) proposed a hybrid optimal control (HOC) framework composed of the inner-loop and the outer-loop to address both continuous and categorical variables, which has been widely used to formalize complex mission planning problems mathematically. In their approach, the outer-loop determines an optimal sequence of targets and the inner-loop solver optimizes the trajectories for the corresponding sequence. A number of studies that applied the HOC framework to design of space missions – such as multiple asteroid missions (Conway et al. 2007; Wall and Conway 2009), interplanetary missions (Englander et al. 2012; Chilan and Conway 2013), and debris removal missions (Yu et al. 2015) – can be found in the literature.

Although the MINLP-based heuristic algorithms and the HOC framework were successfully applied to the multi-target rendezvous problem, exploring optimal rendezvous trajectories between multiple targets is still challenging. It is known that the optimal multi-target rendezvous problem involves extremely large search space (Izzo et al. 2007). To obtain an effective (or, close to optimal) solution of the rendezvous problem within reasonable consumption of computational resource, heuristic algorithms such as the GA and particle swarm optimization (PSO) have been widely used recently. For instance, all the referenced studies employed heuristic algorithms to optimize all design variables for Lambert rendezvous trajectories simultaneously (Zhang et al. 2014; Zhang et al. 2015; Conway et al. 2007; Wall et al. 2009; Englander et al. 2012; Chilan et al. 2013; Yu et al. 2015). It should be noted that the heuristic algorithms are vulnerable to premature convergence and do not guarantee the convergence to the global optimum in general. Chen et al. (2013) pointed out that the heuristic-only approach can fail to find a global solution, even for a single-target rendezvous problem. It is predictable that the effectiveness of the heuristic algorithms will be diminished as the number of targets increases. Therefore, the development of a solution procedure that can smartly handle the multi-target rendezvous problem – including reduction in search space and cost-effective optimum search – is an area that requires attention.

## III. Problem Definition

This section provides the mathematical formulation of the optimal multiple target Lambert rendezvous problem. We first introduce an optimal single-target rendezvous problem ($\mathbf{P}_S$), which is a key component to address the multi-target problem. Then its extension to multi-target rendezvous ($\mathbf{P}_M$) is discussed.

## A. Single-Target Lambert Rendezvous

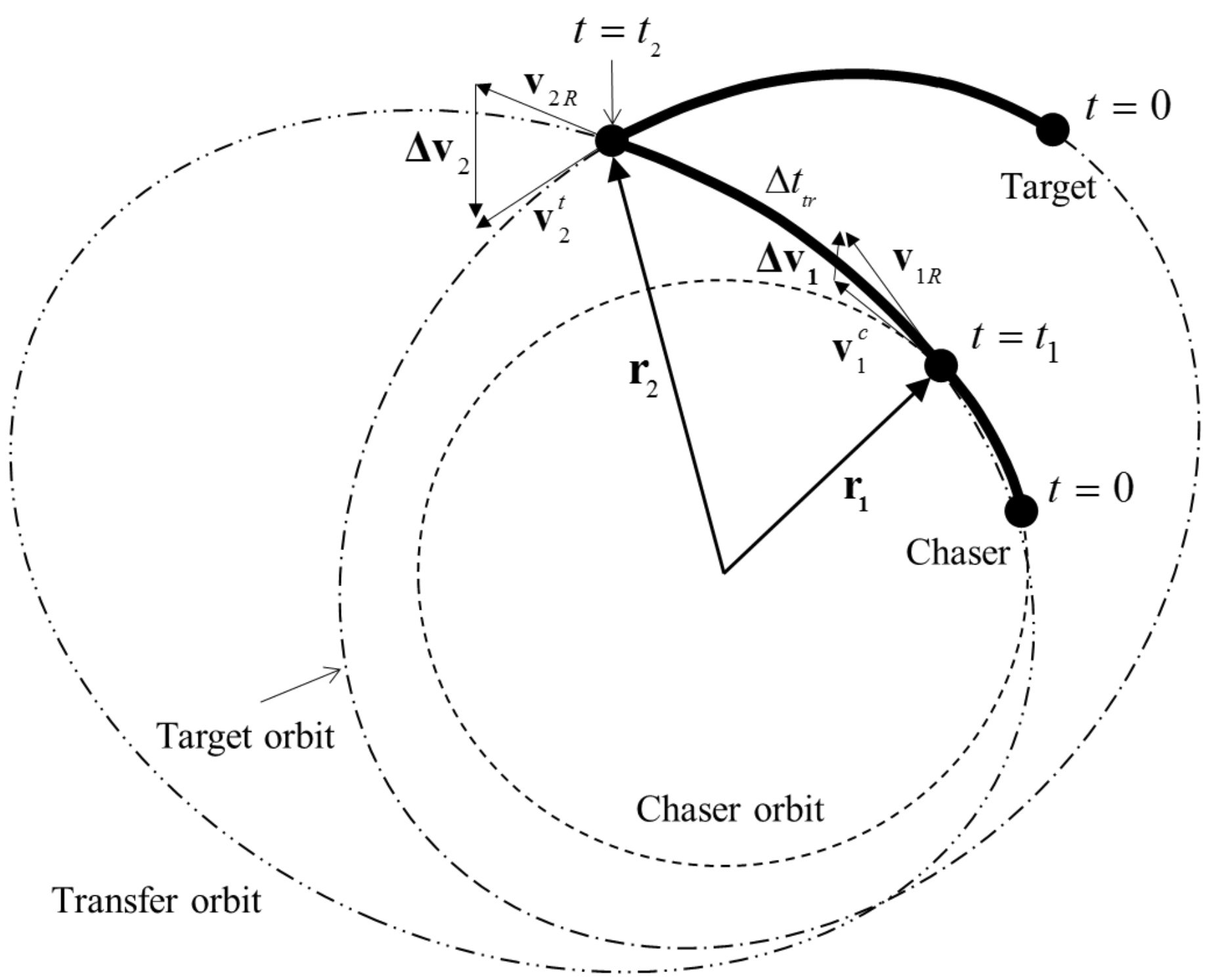

**Figure 1: Single-target Lambert rendezvous problem**

Fig. 1 illustrates a single-target Lambert rendezvous problem (Chen et al. 2013). Initially ($t$ = 0) the chaser and the target are moving in different orbits. After coasting period ($t = t_1$), the first impulsive velocity increment ($\Delta\mathbf{v}_1$) is imposed on the chaser to start an orbital maneuver, which ensures the chaser to meet with the target after transfer time $\Delta t_{tr}$ ($t = t_1 + \Delta t_{tr} = t_2$). The second velocity increment ($\Delta\mathbf{v}_2$) is applied on the chaser spacecraft when the positions of the two objects coincide to make its velocity identical to that of the target (rendezvous condition). Assume that 1) the vectors representing position/velocity of the target are given as functions of time ($\mathbf{r}^t(t)$ / $\mathbf{v}^t(t)$), and 2) the motions of the spacecraft are subject to Newton's law of gravity (Battin 1999). Similarly, given initial position/velocity vectors of the chaser ($\mathbf{r}_0^c, \mathbf{v}_0^c$), its state vectors at $t = t_1$ are expressed as $\mathbf{r}_1^c (= \mathbf{r}_1) = F^c \mathbf{r}_0^c + G^c \mathbf{v}_0^c$ and $\mathbf{v}_1^c = F_t^c \mathbf{r}_0^c + G_t^c \mathbf{v}_0^c$, where $F^c, G^c, F_t^c, G_t^c$ are Lagrange coefficients determined by $\mathbf{r}_0^c, \mathbf{v}_0^c$, and $t_1$.[3]

The total velocity increment necessary to complete the rendezvous can be obtained by solving a multiple-revolution Lambert's problem. The Lambert's problem determines an orbit that has a specific transfer time ($\Delta t_{\text{tr}}$) and initial/final position vectors ($\mathbf{r}_1$ and $\mathbf{r}_2$). The transfer angle – departing from $\mathbf{r}_1$ and arriving at $\mathbf{r}_2$ – can be smaller than 360 deg. (zero-revolution) or greater than 360 deg. (multiple-revolution) depending on the geometry and transfer time of a problem instance. There are a number of published methodologies to solve the multiple-revolution Lambert's problem (Gooding 1990), which provide the required velocity for a spacecraft at the initial point ($\mathbf{v}_{1R}^m = \mathbf{L}_1^m(\mathbf{r}_1, \mathbf{r}_2, \Delta t_{tr})$) and the final position ($\mathbf{v}_{2R}^m = \mathbf{L}_2^m(\mathbf{r}_1, \mathbf{r}_2, \Delta t_{tr})$) for a given number of revolutions ($m$). The optimal single-target rendezvous problem ($\mathbf{P}_\text{S}$) is defined as follows.

[$\mathbf{P}_\text{S}$: Optimal Single-Target Lambert Rendezvous]

$$\min_{t_1, \Delta t_{tr}, m} J_S = (\| \Delta \mathbf{v}_1 \| + \| \Delta \mathbf{v}_2 \|) \tag{1}$$

subject to,

$$\mathbf{r}_1^c = F^c(\mathbf{r}_0^c, \mathbf{v}_0^c, t_1) \cdot \mathbf{r}_0^c + G^c(\mathbf{r}_0^c, \mathbf{v}_0^c, t_1) \cdot \mathbf{v}_0^c \tag{2}$$

$$\mathbf{v}_1^c = F_t^c(\mathbf{r}_0^c, \mathbf{v}_0^c, t_1) \cdot \mathbf{r}_0^c + G_t^c(\mathbf{r}_0^c, \mathbf{v}_0^c, t_1) \cdot \mathbf{v}_0^c \tag{3}$$

$$\mathbf{r}_2^t = \mathbf{r}^t(t_2) = \mathbf{r}^t(t_1 + \Delta t_{tr}) \tag{4}$$

$$\mathbf{v}_2^t = \mathbf{v}^t(t_2) = \mathbf{v}^t(t_1 + \Delta t_{tr}) \tag{5}$$

$$\Delta \mathbf{v}_1 = \mathbf{L}_1^m(\mathbf{r}_1^c, \mathbf{r}_2^t, \Delta t_{tr}) - \mathbf{v}_1^c \tag{6}$$

$$\Delta \mathbf{v}_2 = \mathbf{v}_2^t - \mathbf{L}_2^m(\mathbf{r}_1^c, \mathbf{r}_2^t, \Delta t_{tr}) \tag{7}$$

The objective of the problem presented in Eq. (1) is to minimize sum of velocity increments

---

[3] Throughout the paper, subscripts 1 and 2 represent the beginning and end of the transfer maneuver, and subscripts $c$ and $t$ denote the chaser and target.

(for departure and arrival) to conduct the rendezvous with the target. Design variables are the initial time to start the rendezvous ($t_1$), the transfer time ($\Delta t_{tr}$), and the number of revolutions of the rendezvous trajectory ($m$). Eqs. (2)-(3) express the velocity and position of the chaser spacecraft at the beginning of the orbital transfer (at $t = t_1$) and Eqs. (4)-(5) represent the states of the target at the end of the transfer ($t = t_2 = t_1 + \Delta t_{tr}$)). The expression for the departure and arrival velocity increments are presented in Eqs. (6)-(7).

## B. Multi-Target Lambert Rendezvous

The optimal multi-target Lambert rendezvous problem is formulated by extending the single-target problem (**P**$_S$) expressed in Eqs. (1)-(7). The problem determines the optimal rendezvous sequence and associated trajectories for the series of single-target Lambert rendezvous with given targets. Suppose that the position/velocity of the chaser spacecraft at the beginning of the mission are $\mathbf{r}_0^c / \mathbf{v}_0^c$ and those of the target $i$ ($i$ = 1, … , $N$; $N$ is the number of targets) at time $t$ are expressed as $\mathbf{r}^{t,i}(t) / \mathbf{v}^{t,i}(t)$. The optimal multi-target Lambert rendezvous problem (**P**$_M$) is defined as follows.

[**P**$_M$: Optimal Multi-Target Lambert Rendezvous]

$$\min_{\mathbf{x},\,\mathbf{q}} J_M = \sum_{i=1}^{N} (\| \Delta \mathbf{v}_{1,i} \| + \| \Delta \mathbf{v}_{2,i} \|) \tag{8}$$

subject to,

$$\mathbf{x} = [\boldsymbol{\tau}_1, \Delta\boldsymbol{\tau}_{tr}, \boldsymbol{\mu}] \tag{9}$$

$$\boldsymbol{\tau}_1 = [t_{1,1}, \cdots, t_{1,N}] \tag{10}$$

$$\Delta\boldsymbol{\tau}_{tr} = [\Delta t_{tr,1}, \cdots, \Delta t_{tr,N}] \tag{11}$$

$$\boldsymbol{\mu} = [m_1, \cdots, m_N] \tag{12}$$

$$\mathbf{q} = [q_1, \cdots, q_N] \tag{13}$$

$$t_{0,i} = t_{2,(i-1)} + t_{ser} = t_{1,(i-1)} + \Delta t_{tr,(i-1)} + t_{ser} \quad (i \ge 2), \;\; t_{0,1} = 0 \tag{14}$$

$$\mathbf{r}_{0,i}^c = \mathbf{r}^{t,q_{(i-1)}}(t_{0,i}) \;\; (i \ge 2), \;\; \mathbf{r}_{0,1}^c = \mathbf{r}_0^c \tag{15}$$

$$\mathbf{v}_{0,i}^{c} = \mathbf{v}^{t,q_{(i-1)}}(t_{0,i}) \;\; (i \geq 2), \;\; \mathbf{v}_{0,1}^{c} = \mathbf{v}_{0}^{c} \tag{16}$$

$$\mathbf{r}_{1,i}^{c} = F^{c}(\mathbf{r}_{0,i}^{c}, \mathbf{v}_{0,i}^{c}, t_{1,i} - t_{0,i}) \cdot \mathbf{r}_{0,i}^{c} + G^{c}(\mathbf{r}_{0,i}^{c}, \mathbf{v}_{0,i}^{c}, t_{1,i} - t_{0,i}) \cdot \mathbf{v}_{0,i}^{c} \tag{17}$$

$$\mathbf{v}_{1,i}^{c} = F_{t}^{c}(\mathbf{r}_{0,i}^{c}, \mathbf{v}_{0,i}^{c}, t_{1,i} - t_{0,i}) \cdot \mathbf{r}_{0,i}^{c} + G_{t}^{c}(\mathbf{r}_{0,i}^{c}, \mathbf{v}_{0,i}^{c}, t_{1,i} - t_{0,i}) \cdot \mathbf{v}_{0,i}^{c} \tag{18}$$

$$\mathbf{r}_{2,i}^{t} = \mathbf{r}^{t,q_i}(t_{2,i}) = \mathbf{r}^{t,q_i}(t_{1,i} + \Delta t_{tr,i}) \tag{19}$$

$$\mathbf{v}_{2,i}^{t} = \mathbf{v}^{t,q_i}(t_{2,i}) = \mathbf{v}^{t,q_i}(t_{1,i} + \Delta t_{tr,i}) \tag{20}$$

$$\Delta\mathbf{v}_{1,i} = \mathbf{L}_{1}^{m_i}(\mathbf{r}_{1,i}^{c}, \mathbf{r}_{2,i}^{t}, \Delta t_{tr,i}) - \mathbf{v}_{1,i}^{c} \tag{21}$$

$$\Delta\mathbf{v}_{2,i} = \mathbf{v}_{2,i}^{t} - \mathbf{L}_{2}^{m_i}(\mathbf{r}_{1,i}^{c}, \mathbf{r}_{2,i}^{t}, \Delta t_{tr,i}) \tag{22}$$

The objective of this problem is to minimize the sum of ΔV values used to complete all ($N$) rendezvouses, which is presented in Eq. (8). The design variables of the problem are presented in Eqs. (9)-(13). Index $i$ represents the order of rendezvous task and variables $t_{1,i}$, $\Delta t_{tr,i}$, $m_i$, and $q_i$ denote departure time, transfer time, number of revolution, and target index associated with the $i^{th}$ rendezvous, respectively. Initial time, position, and velocity for the $i^{th}$ rendezvous task ($t_{0,i}, \mathbf{r}_{0,i}^{c}, \mathbf{v}_{0,i}^{c}$) are determined by the final time of the previous ($(i\text{-}1)^{th}$) rendezvous and the required service time ($t_{ser}$) at the target, which are recursively defined in Eqs. (14)-(16). The states of the chaser spacecraft at the beginning of the orbital transfer, the states of the target at the end of the transfer, and the expressions for the departure and arrival velocity increments for the $i^{th}$ rendezvous are presented in Eqs. (17)-(22), which are the "multi-target version" counterparts of Eqs. (2)-(7). Fig. 2 illustrates the multi-target Lambert rendezvous process for the case of $N$ = 4 and $\mathbf{q}$ = (2, 1, 3, 4).

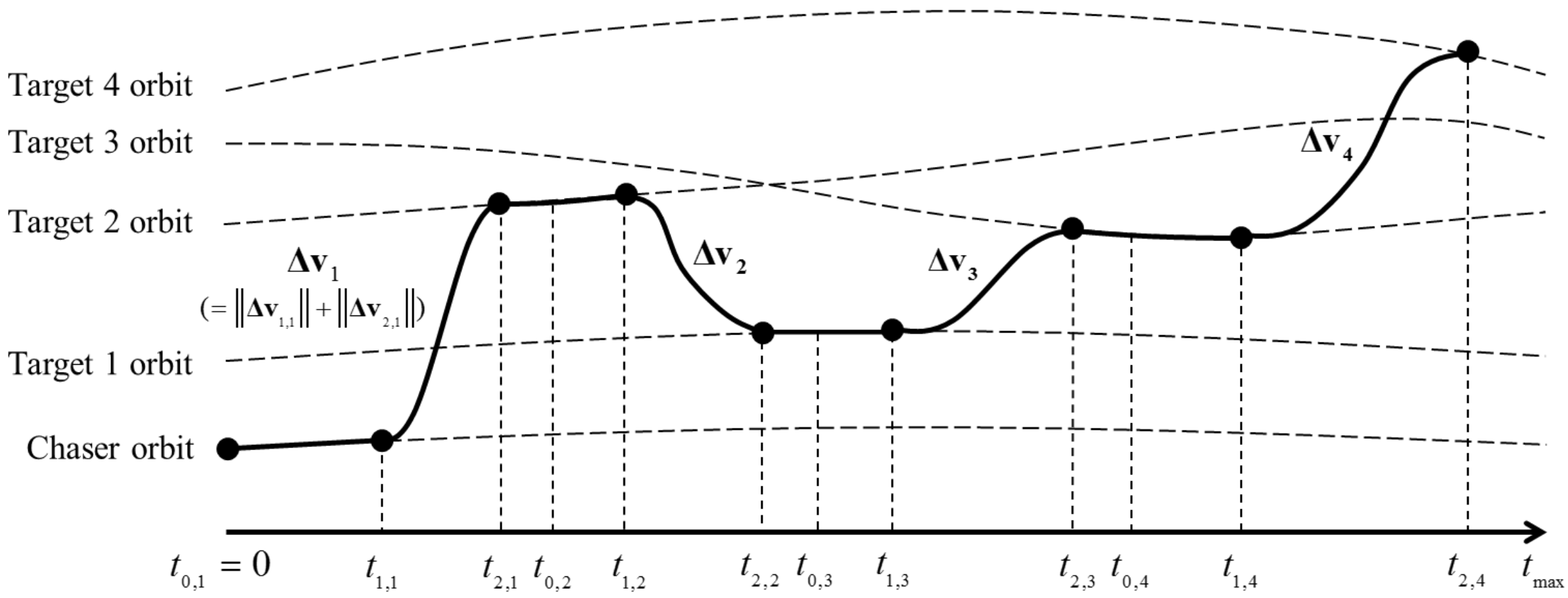

**Figure 2: Example of multi-target Lambert rendezvous process (N = 4, $\mathbf{q}=(2,1,3,4)$ )**

# IV. Framework to Solve an Optimal Multi-Target Lambert Rendezvous

As was discussed in the previous section, the optimal multi-target Lambert rendezvous problem ($\mathbf{P}_M$) can be classified as MINLP, which is very difficult to handle directly. This section proposes a framework for solving the problem composed of two distinct phases. The structure of the proposed framework is outlined in Fig. 3.

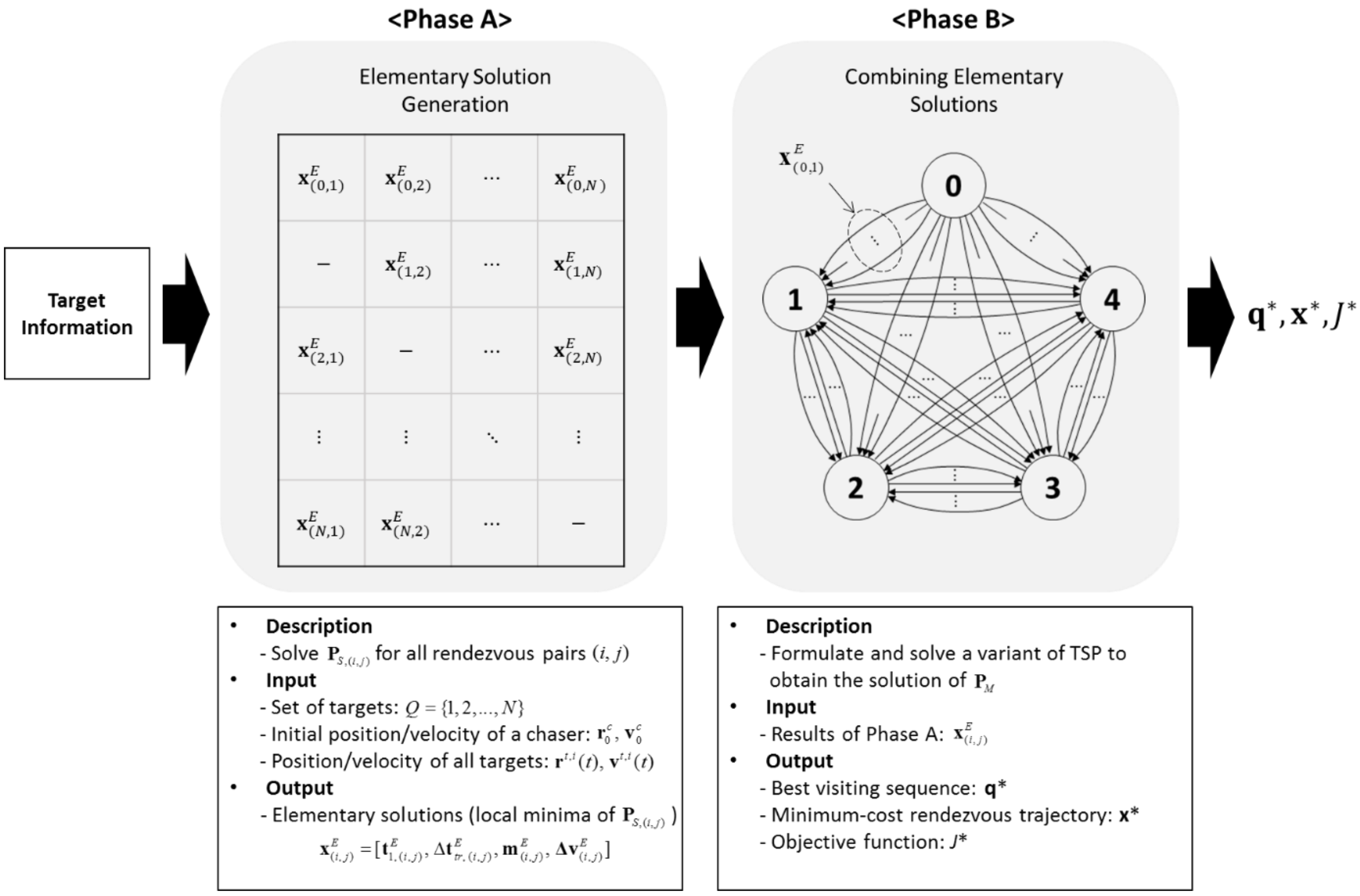

**Figure 3: Framework to solve the multi-target Lambert rendezvous problem ($\mathbf{P}_M$)**

In the first phase of the framework, a series of single target rendezvous problems ($\mathbf{P}_S$) for every departure/arrival object pair are solved to obtain all of their local minima, which are used as the elementary solutions for the original multi-target rendezvous problem. Each elementary solution represents a candidate rendezvous trajectory between the object pair that has its own cost, transfer time, and time window. The second phase seeks for the best rendezvous sequence ($\mathbf{q}^*$) and trajectories ($\mathbf{x}^*$) of the original problem using the elementary solutions obtained in the first phase. A new variant of TSP that considers multiple arcs (corresponding to the elementary solutions) and arc time window constraints is introduced to find the best combination of the elementary solutions. Note that all computations for the rendezvous trajectory optimization, which is the NLP part of the problem, are carried out in the first phase, while the second phase conducts the optimal sequencing, the ILP part of the problem. The following subsections provide the details on the steps of the proposed framework.

### A. Phase A: Elementary Solution Generation

Phase A obtains all local minima of every single-target rendezvous problem ($\mathbf{P}_S$) instantiated by specified departure/arrival objects, which provide the elementary solutions of the overall problem (Bang and Ahn 2016)[4]. We first explore the patterns in the locations of local minima of $\mathbf{P}_S$.

Fig. 4 presents a sample contour plot of cost function for a single-target Lambert rendezvous (the sum of departure and arrival ΔV's) in $t_1 - \Delta t_{tr}$ plane. The orbits of the chaser and the target are elliptic and non-coplanar. While the contours are diverse in their shape depending on orbital elements of the chaser and target (e.g. semi-major axis, eccentricity, and inclination angle), one can observe some patterns common to the locations of local minima. Frist, the cost function is extremely high along equally spaced straight lines (or, "walls") of two different types (type A and type B). Secondly, there are multiple minima located in regions separated by the walls.

The locations of walls (illustrated as lines A, A', B and B' in Fig. 4) depend on the transfer geometry associated with the initial and final position vectors ($\mathbf{r}_1$ and $\mathbf{r}_2$) of the rendezvous maneuver. Fig. 5 presents the orbital transfer geometries that can cause large velocity increments. In the figure, $l$ is the line of intersection between the orbital planes and $\theta_1 / \theta_2$ are the angles between $l$ and $\mathbf{r}_1 / \mathbf{r}_2$, respectively. Fig. 5-(a) (zero transfer angle for a coplanar rendezvous) visualizes the geometry corresponding to type A walls (lines A and A' in Fig. 4). A rendezvous maneuver for this geometry requires very large velocity increment since the tangential velocity of the chaser should be lost at the beginning and then recovered at the end of the rendezvous. Therefore, the local minima are not located on or near the type A walls. The horizontal and vertical distances between two adjacent type A walls are given as

$$\delta t_{1,A} = \frac{2\pi}{n_1 - n_2} \tag{23}$$

$$\delta(\Delta t_{tr,A}) = \frac{2\pi}{n_2} \tag{24}$$

where $n_1$ and $n_2$ are mean motions of the chaser and the target.

---

[4] A preliminary version of the procedure to solve the single-target rendezvous optimization problem introduced in this subsection was presented in authors' previous conference paper (Bang and Ahn 2016).

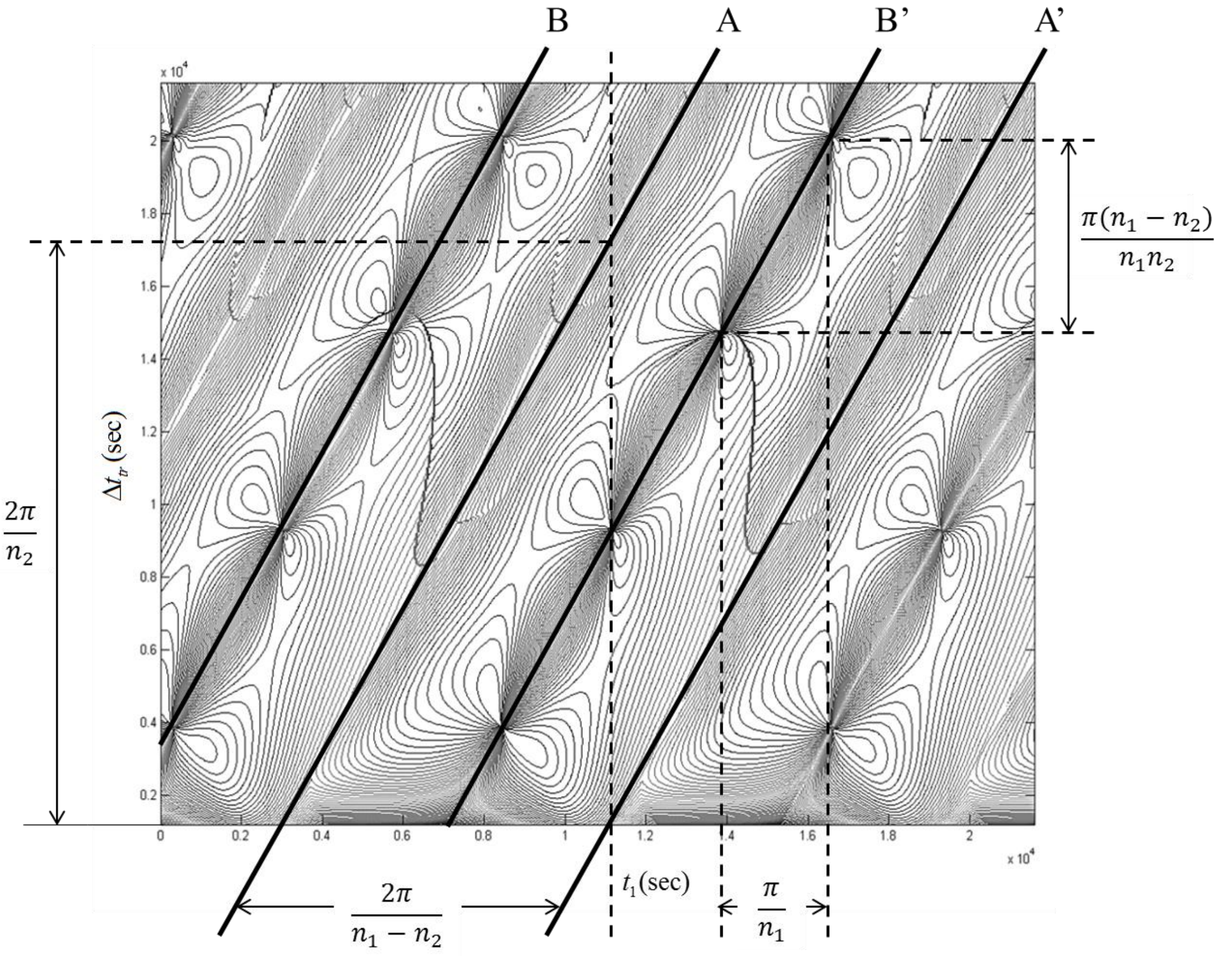

**Figure 4: Contour plot of cost function for a single-target Lambert rendezvous**

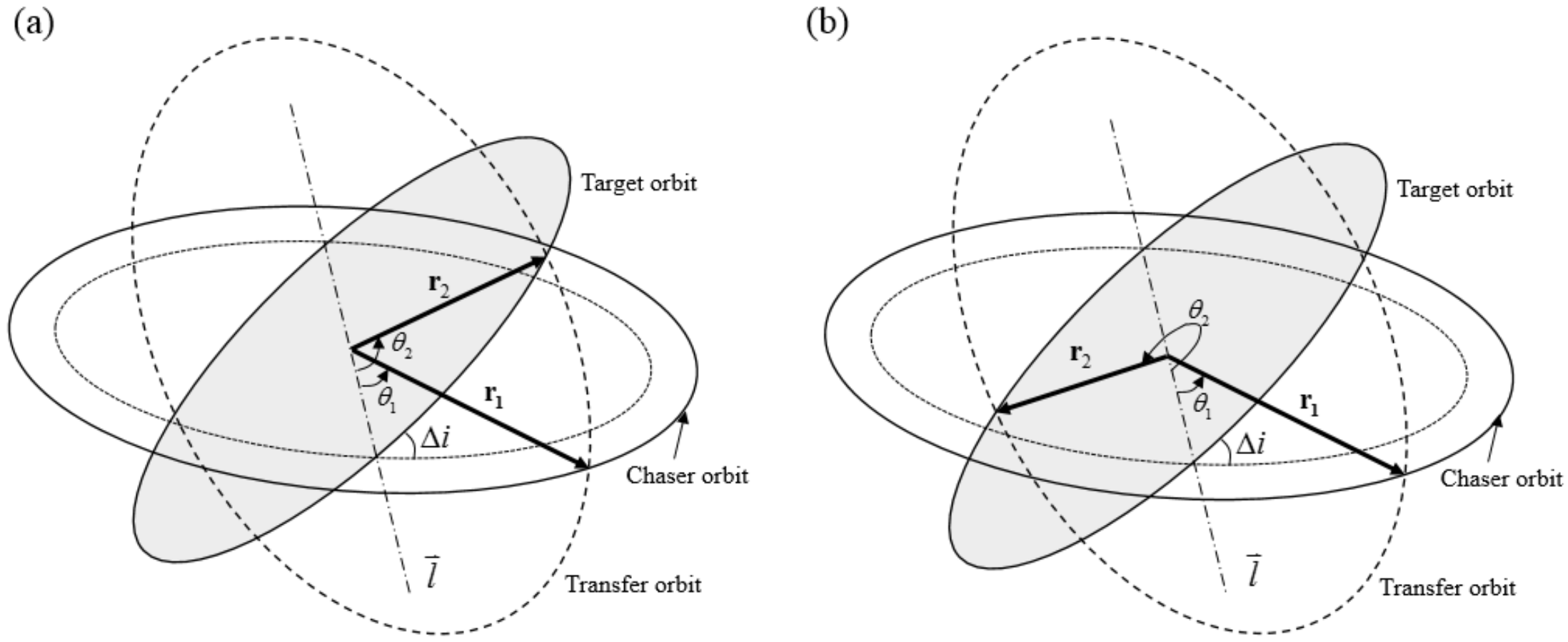

**Figure 5: Orbital transfer geometries associated with type A and type B walls**

Type B walls appear only in the contour plots for non-coplanar rendezvous, which is presented in Fig. 5-(b). The value of cost function changes periodically along the lines B and B' of Fig. 4. The required ΔV is minimized when $(\theta_1, \theta_2)$ equals (0°, 180°) or (180°, 0°), where $\mathbf{r}_1$ and $\mathbf{r}_2$ lie on the intersection of the two orbital planes. The local minima are placed on or near the type B walls with horizontal and vertical spaces as

$$\delta t_{1,B} = \frac{\pi}{n_1} \tag{25}$$

$$\delta(\Delta t_{tr,B}) = \frac{\pi(n_1 - n_2)}{n_1 n_2} \tag{26}$$

The procedure to find all minima of a single-target Lambert rendezvous problem is developed based on the aforementioned characteristics of the contour plot. A solution space partitioning technique is used to narrow down the exploration region and a gradient-based algorithm is implemented to obtain the minima. The solution space partitioning splits the current exploration space into multiple subspaces (exploration spaces of the next iteration step) until the subspace has at most one local minimum. Note that, when we explore the local minima of a subspace, constraints on the ranges of $t_1$ and $\Delta t_{tr}$ corresponding to the definition of the subspace are added to problem $\mathbf{P}_S$ as follows:

$$t_{1,\min} \le t_1 \le t_{1,\max} \tag{27}$$

$$\Delta t_{tr,\min} \le \Delta t_{tr} \le \Delta t_{tr,\max} \tag{28}$$

A gradient-based algorithm with multiple initial points, which leverages the characteristics of the contour, is introduced to solve the single-target Lambert rendezvous problem with additional constraints ($\mathbf{P}_S$ with Eqs. (27)-(28)). When the procedure is conducted with initial guess at four vertices of an arbitrary subspace, the results can be categorized into the following three cases, which are presented in Fig. 6.

Case (a): If there is no local minimum of the single-target rendezvous problem ($\mathbf{P}_S$) inside the subspace, all the solution(s) of the constrained problem ($\mathbf{P}_S$ with additional subspace constraint) are located on the subspace boundary.

Case (b): If there is only one local minimum of $\mathbf{P}_S$ and there is no wall inside the subspace, all the solutions of the constrained problem will converge on the local minimum.

Case (c): If there are two or more local minima of $\mathbf{P}_S$ or there exists a wall inside the subspace, some solutions converge to the minima and others are located on the boundary.[5]

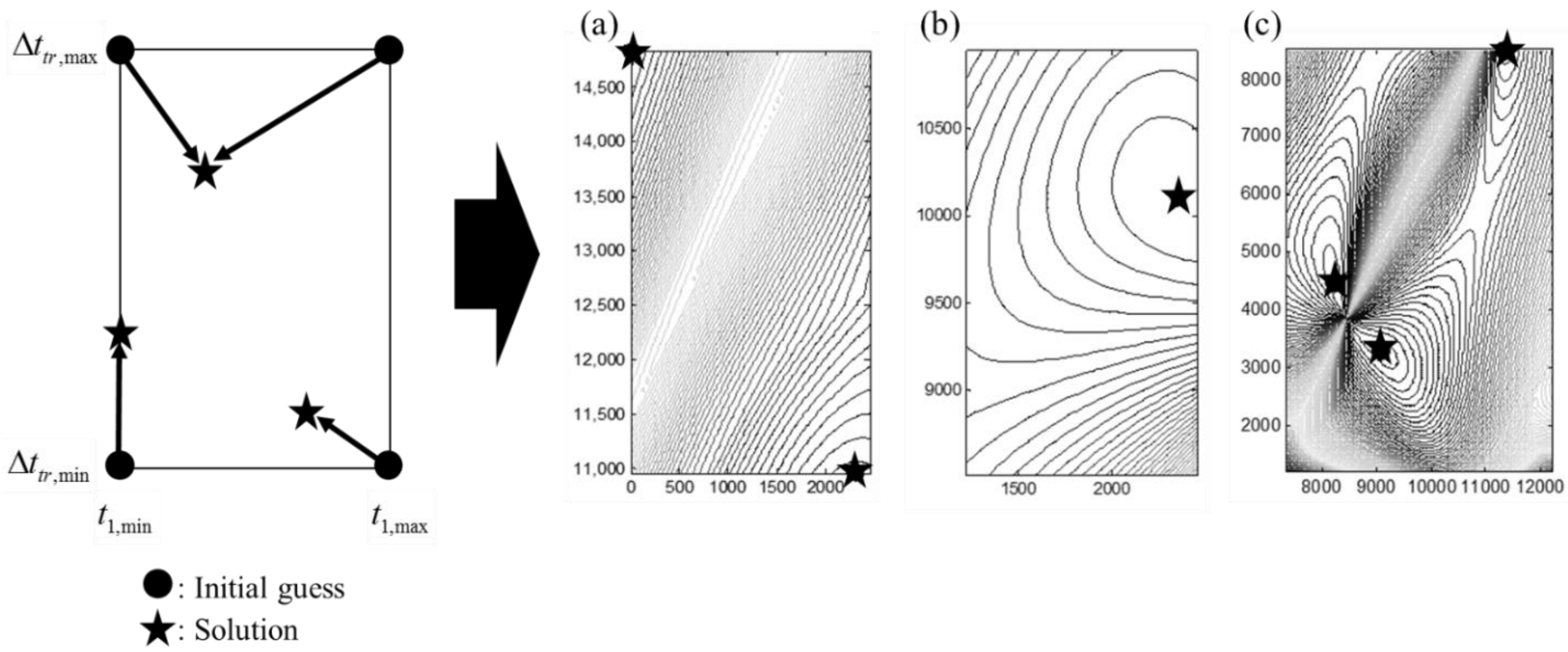

**Figure 6: Gradient-based algorithm implemented at four grid vertices**

Cases (a) and (b) leads to the conclusions that there is no or just one local minimum in the subspace and additional split is not necessary. For Case (c), on the other hand, not all the local minima have been identified and further exploration of split solution space is required. The procedure to find

---

[5] We set the initial grids to define the subspaces sufficiently fine and can assume that the number of type A wall inside of the subspace is no larger than one.

all local minima of a single-target rendezvous problem is summarized in Table 1.

**Table 1: Procedure to find all local minima of $\mathbf{P}_S$**

| Step | Task Description |
|---|---|
| Step 1 | Partition the whole solution space into subspaces using initial grids. Note that the size of each subspace should be smaller than $(\delta t_{1,A} \times \delta(\Delta t_{tr,A}))$. |
| Step 2 | For each grid, implement the gradient-based algorithm (e.g. sequential quadratic programming (SQP)) with initial guesses of four vertex points. The lower and upper bounds for design variables are set to be identical to the boundary of the grid. |
| Step 3 | If all the solutions are located on the boundary (Case (a)), stop exploring the subspace (eliminate the subspace) since there is no local minimum within it. Otherwise, go to Step 4. |
| Step 4 | If all the solutions converge on a single point inside the boundary (Case (b)), save the point as a local minimum (elementary solution) and stop exploring the subspace (eliminate the subspace). Otherwise, go to the Step 5. |
| Step 5 | Split the exploring subspace into smaller pieces because the identification of all the local minima is not guaranteed (Case (c)). Repeat steps 2 to 4 for new subspaces until every subspace is eliminated |

## B. Phase B: Solving Multi-Target Problem with Elementary Solutions

Phase B of the proposed framework solves the multi-target rendezvous problem ($\mathbf{P}_M$) by combining elementary solutions prepared in Phase A. The first step is to construct a graph composed of 1) nodes representing the objects, 2) arcs representing the rendezvous trajectories between two objects, and 3) costs associated with the arcs. Let $G = (V, A)$ be a directed graph where $V$ $(= \{0\} \cup Q = \{0, 1, \ldots, N\})$ is the set of nodes (0: initial chaser orbit, $Q = \{1, \ldots, N\}$: targets) and $A$ $(= \bigcup_{(i,j)\in V\times V,\, i\neq j} A_{(i,j)})$ is a set of arcs. Note that there is an infinite number of possible rendezvous trajectories with different costs, departure times, and transfer times between each object pair; hence, it is not appropriate to confine the path between each pair of nodes to a single arc. Instead, the elementary solutions $\mathbf{x}^E_{(i,j)}$, which represent a set of candidate rendezvous trajectories from object $i$ to object $j$, are used to define a set of multiple arcs from node $i$ to node $j$ as

$$A_{(i,j)} = \{(i,\ j)^p \mid 1 \le p \le \left|\mathbf{x}^E_{(i,j)}\right|\} \tag{29}$$

Fig. 7 illustrates an example of a graph with four targets (*N*=4) constructed using the elementary solutions. In the graph, the proposed multi-target rendezvous problem ($\mathbf{P}_M$) is equivalent to the process to determine the path in the graph that starts at node 0 and visits all other nodes with one of multiple arcs while minimizing the sum of costs (ΔV values).

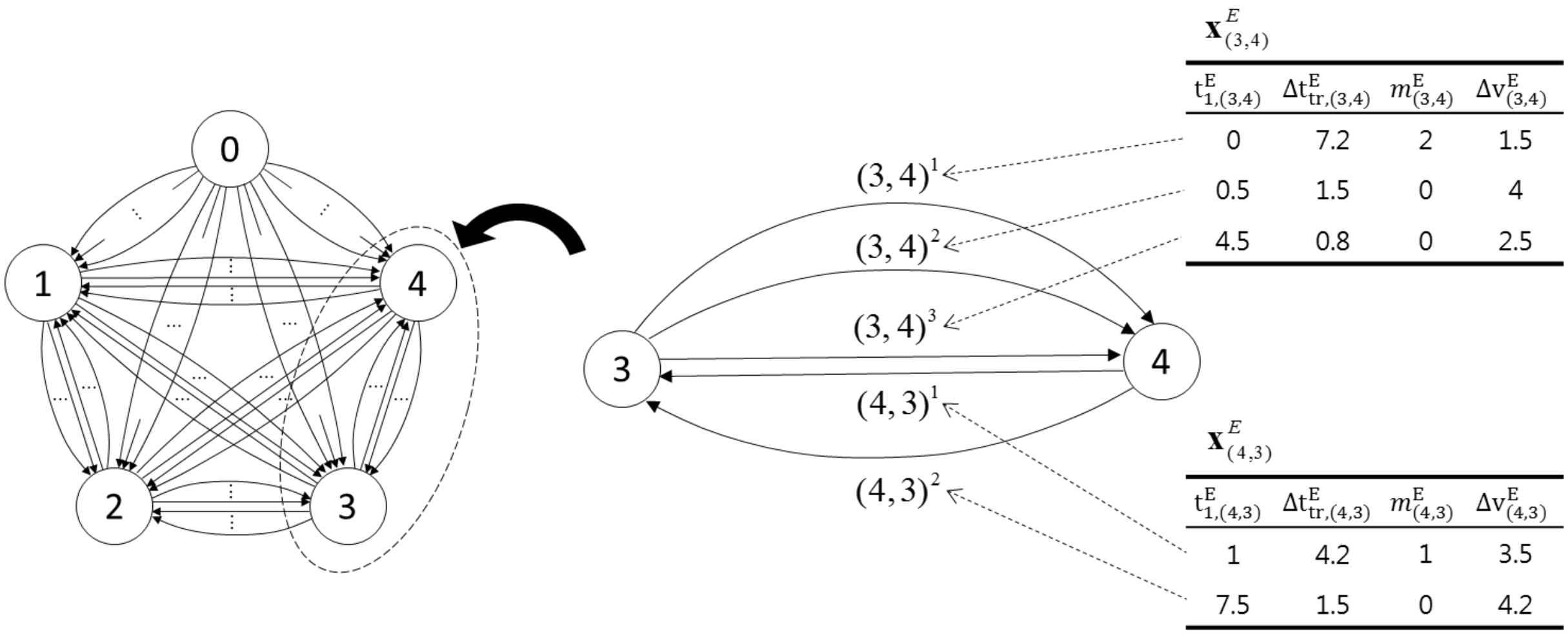

**Figure 7: Graph for multi-target rendezvous problem ($\mathbf{P}_M$) – an illustrative example (*N* = 4)**

Determination of an optimal path in the proposed graph can be interpreted as a variant of TSP formulated by considering two additional features: 1) there are multiple arcs between two nodes (i.e. departure and arrival objects), and 2) an arc is available only during a specific time interval referred to as *arc time window*. Note that each elementary solution (arc in the graph) has its own departure time ($t_1$). Let $t$ denote the time to finish the service at a specific node. If the service at a node finishes before the departure time associated with an arc originating from the node ($t \in [0, t_1]$), the spacecraft can wait and depart to rendezvous with the next target using the arc. Otherwise ($t > t_1$), the spacecraft cannot use the arc (rendezvous opportunity missed). Therefore, the time window for the arc is defined as $[0, t_1]$.

This situation is formulated as a new variant of TSP considering multiple arcs and arc time windows. The decision variable $y_{(i,j)^p}$ (binary) takes value of 1 if the spacecraft travels from node $i$ to node $j$ using arc $(i, j)^p$ – the $p^{th}$ arc associated with the departure/arrival objects, and 0 otherwise. In addition, $c_{(i,j)^p}$, $t_{1(i,j)^p}$, and $\Delta t_{tr(i,j)^p}$ denote the cost, the departure time (latest available time), and the transfer time associated with $(i, j)^p$, respectively. The mathematical formulation for Phase B of the proposed framework is presented as follows:

[$\mathbf{P}_{\mathrm{ME}}$: Optimal Multi-Target Rendezvous with Elementary Solutions]

$$\min_{\mathbf{y}} J_{ME} = \sum_{i \in V} \sum_{j \in V} \sum_{p=1}^{|A_{(i,j)}|} c_{(i,j)^p} y_{(i,j)^p} \tag{30}$$

subject to,

$$\sum_{i \in V} \sum_{p=1}^{|A_{(i,j)}|} y_{(i,j)^p} = 1, \forall j \in V \tag{31}$$

$$\sum_{j \in V} \sum_{p=1}^{|A_{(i,j)}|} y_{(i,j)^p} = 1, \forall i \in V \tag{32}$$

$$\sum_{j \in V} \sum_{p=1}^{|A_{(j,i)}|} (t_{1(j,i)^p} + \Delta t_{tr(j,i)^p}) y_{(j,i)^p} + t_{ser} \le \sum_{j \in V} \sum_{p=1}^{|A_{(i,j)}|} t_{1(i,j)^p} y_{(i,j)^p}, \forall i \in Q \tag{33}$$

$$y_{(i,j)^p} \in \{0,1\}, \forall i, j \in V, \ \forall p \in \{1,\ldots,|A_{(i,j)}|\} \tag{34}$$

where $\mathbf{y}( = [\cdots, y_{(i,j)^p}, \cdots]^T)$ in Eq. (30) is the vector collecting the decision variables.

The objective function presented in Eq. (30) is to minimize sum of costs associated with the path. Eqs. (31)-(32) represent the constraints that each node is visited exactly once and there is exactly one departure from each node. Eq. (33) guarantees that the completion of service at the node should occur earlier than the latest available time of the following arc. Moreover, Eq. (33) prevents sub-tours. Finally, the integrity constraint is presented in Eq. (34). Note that $\mathbf{P}_{\mathrm{ME}}$ is a routing problem that considers both of multiple arcs assigned between two nodes and the arc time windows (*Traveling Salesman Problem with Multiple Arcs and Arc Time Windows*), which has not yet been addressed in operations research field as far as the authors know.

# V. Case Study

This section presents the case studies demonstrating the effectiveness of the proposed framework to solve an optimal multi-target Lambert rendezvous problem. Subsection V-A validates the algorithm to find a set of local minima for a single-target Lambert rendezvous (Phase A). Subsection V-B solves multi-target problems by adopting the whole framework introduced in this paper (Phases A and B). All the numerical results are compared with the solutions obtained by the generic algorithm.

## A. Case A: Validation for Single-Target Lambert Rendezvous

The proposed algorithm based on *subspace partitioning* and *gradient search* to find minima of a single target Lambert rendezvous was applied to a test problem presented in Chen et al. (2013). Table 2 summarizes the orbital elements of chaser/target and the constraint for the problem. The size of grids defining the initial subspace is set as: $(0.9\times\delta t_{1,A})\times(0.9\times\delta(\Delta t_{tr,A}))$. A subspace subject to additional split (Case (c) of Section IV-A) is divided into four by halving its width/height. The procedure was implemented in MATLAB and the fmincon function of MATLAB optimization toolbox was adopted as a gradient-based optimizer for local minimum search.

**Table 2: Orbit elements and parameters for Case A (Chen et al. 2013)**

| **Orbital elements** | **Chaser** | **Target** |
|---|---:|---:|
| Semi-major axis, *km* | $R_E$+2,000 | $R_E$+36,000 |
| Eccentricity, - | 0.002 | 0.0002 |
| Inclination, *deg* | 60.00 | 55.00 |
| Right ascension of the ascending node, *deg* | 30.00 | 35.00 |
| Argument of perigee, *deg* | 0.00 | -20.00 |
| True anomaly, *deg* | 0.00 | 30.00 |
| **Parameter** | | **Value** |
| $t_{max}$, *hr* | | 24 |

Fig. 8 illustrates the test case result. Total 37 local minima (depicted as circles) – including those located on the boundary representing the terminal time constraint – were found through the procedure introduced in Section IV-A. The figure shows that the procedure finds all local minima without any misses. In addition, the global optimum (depicted as a black square) obtained by the procedure is identical to the reported value. Chen et al. (2013) pointed out that the global optimum obtained by the GA for this problem was inconsistent (one of points depicted as crosses). This observation indicates that 1) a single run of GA does not guarantee the global optimum, and 2) a number of optimization runs are required to find the global optimum with acceptably high probability. On the contrary, the proposed approach can find the set of all local minima, which include the global optimum.

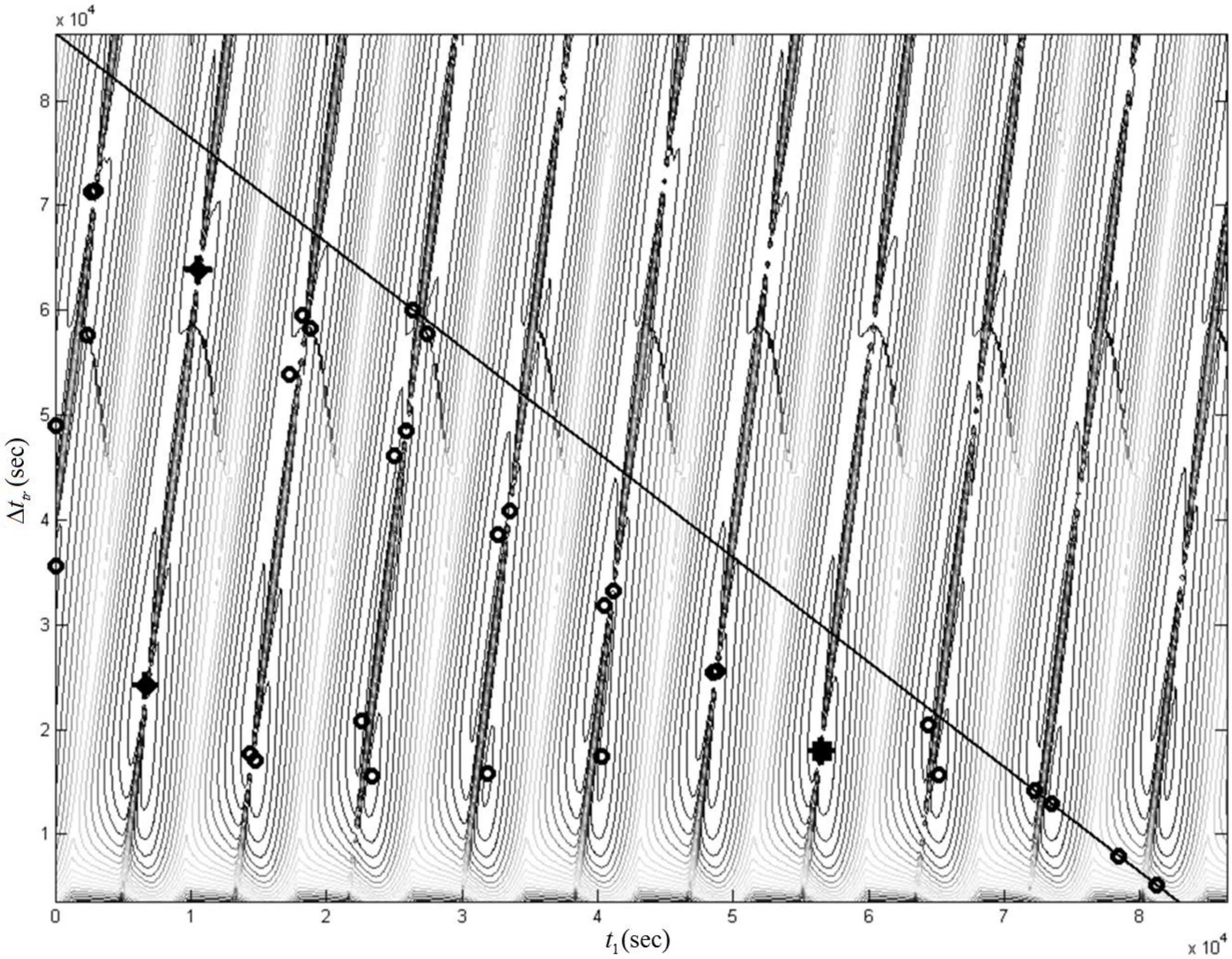

**Figure 8: Results of Case A – local minima and global optimum obtained by the propose procedure and the GA**

## B. Case B: Validation for Multi-Target Lambert Rendezvous

The overall framework proposed in the paper is validated using the multi-asteroid rendezvous case presented in Zhang et al. (2015). The objective of the case is to find the best (min- $\Delta V$ ) rendezvous order and trajectories to visit all asteroid within a specified mission time. Three problem instances with different number of asteroids – 4 (Case B-1), 8 (Case B-2), and 16 (Case B-3) – were considered and the framework introduced in section IV (Phases A and B) is used to solve them. The orbital elements of asteroids and the constraint on the final time for the problem were summarized in Table 3 (Zhang et al. 2015). The TSP with multiple arcs and arc time windows formulated in Phase B was implemented in MATLAB with GUROBI 7.5 solver.

**Table 3: Asteroid orbit elements and mission parameters for Case B (Zhang et al. 2015)**

| ID | Name | Epoch, JD | a, AU | e, - | i, deg | Ω, deg | ω, deg | M, deg |
|---|---|---|---|---|---|---|---|---|
| 1 | 2001-GP2 | 2456600.5 | 1.0377497 | 0.0740190 | 1.27980 | 196.84658 | 111.32136 | 220.55456 |
| 2 | 2007-UN12 | 2456600.5 | 1.0537339 | 0.0604714 | 0.23523 | 216.10320 | 134.34440 | 242.72351 |
| 3 | 2006-JY26 | 2456600.5 | 1.0099643 | 0.0830940 | 1.43929 | 43.48569 | 273.55713 | 57.93400 |
| 4 | 2010-JR34 | 2456600.5 | 0.9593424 | 0.1448353 | 0.72205 | 36.86737 | 316.29050 | 146.08995 |
| 5 | 2009-BD | 2456600.5 | 1.0617043 | 0.0515307 | 1.26705 | 253.32360 | 316.76764 | 115.45581 |
| 6 | 2008-JL24 | 2456600.5 | 1.0382537 | 0.1065636 | 0.55071 | 225.80291 | 281.99557 | 136.76768 |
| 7 | 2008-UA202 | 2456600.5 | 1.0332190 | 0.0684246 | 0.26337 | 21.04400 | 300.95352 | 345.64862 |
| 8 | 2006-BZ147 | 2456600.5 | 1.0235198 | 0.0985837 | 1.40948 | 139.83156 | 94.80897 | 84.23810 |
| 9 | 2009-BK2 | 2456600.5 | 1.0125947 | 0.2128808 | 3.57355 | 126.43126 | 121.32795 | 144.17500 |
| 10 | 2001-CC21 | 2456600.5 | 1.0324714 | 0.2194123 | 4.80881 | 75.58951 | 179.33863 | 326.38628 |
| 11 | 2001-QJ142 | 2456600.5 | 1.0621126 | 0.0862374 | 3.10322 | 184.40414 | 63.93151 | 125.58762 |
| 12 | 2009-OS5 | 2457000.5 | 1.1441219 | 0.0967254 | 1.69487 | 145.37326 | 120.83594 | 173.44777 |
| 13 | 1999-AO10 | 2456600.5 | 0.9115237 | 0.1109538 | 2.62074 | 313.27180 | 7.66737 | 147.50501 |
| 14 | 2013-BS45 | 2457000.5 | 0.9936781 | 0.0838748 | 0.77337 | 83.55062 | 149.70642 | 219.48084 |
| 15 | 2013-NX | 2456800.5 | 1.0323256 | 0.1698532 | 6.32151 | 112.73144 | 312.66316 | 177.28424 |
| 16 | 2012-UV136 | 2457000.5 | 1.0075469 | 0.1389606 | 2.20908 | 210.52858 | 289.49134 | 305.08581 |

| **Mission Parameter** | **Value** |
|---|---|
| $t_0$, JD | 2457023.5 |
| $t_{ser}$, days | 7 |

| **Problem Instance** | **No. of Asteroids** | **Given Target Set** | $t_{max}$**, yr** |
|---|---|---|---|
| Case B-1 | 4 | {1, 2, 3, 4} | 22 |
| Case B-2 | 8 | {1, 2, …, 8} | 34 |
| Case B-3 | 16 | {1, 2, …, 16} | 58 |

Table 4 comparatively exhibits the results of the three test cases obtained using the proposed framework and the mixed-code genetic algorithm (MCGA) with search enhancement presented in Zhang et al. (2015). The optimal visiting orders for Case B-1 (rendezvous with 4 asteroids) obtained by two different methods were identical and their optimal objective function values show no significant difference. For Case B-2, however, the solutions found by different methodologies were very different – both in terms of rendezvous sequence and the best objective function ([6, 8, 5, 1, 4, 3, 7, 2] and 16.406 *km*/*s* for the *proposed framework* and [8, 5, 1, 4, 7, 3, 2, 6] / 19.153 *km*/*s* for the *MCGA with search enhancement*). The improvement in the objective function achieved by the proposed framework was 14.3 %. In Case B-3 (rendezvous with 16 asteroids), the optimal solution found by the proposed framework outperformed the reported result obtained by the MCGA with search enhancement – more significantly than previous cases. The visiting sequences obtained by two methodologies were totally different (proposed framework: [12, 11, 1, 5, 13, 3, 7, 2, 4, 15, 16, 6, 14, 8, 9, 10], MCGA: [11, 9, 1, 4, 12, 10, 3, 7, 2, 5, 6, 16, 14, 8, 15, 13]), and the reduction in objective function was 35.6 % (37.325 *km*/*s* versus 58.977 *km*/*s*). This comparison result supports the aforementioned limitation of existing heuristics based approach in its effectiveness for rendezvous with many targets. The number of Lambert routine calls for the two methods, which provide the proxy for the resource consumption, are compared as well. The results indicate that the computational resource spent by the proposed method is about 3-4 times larger than that for the MCGA.

The performance enhancement of the proposed framework will be one of interesting follow-on study subjects. The improvement can be made by both revising the solution procedure and accelerating the Lambert routine, which is the bottleneck process that consumes the majority of computation time (Ahn and Lee 2013; Ahn et al. 2015; Lee et al. 2016).

**Table 4: Objective function and computational load for Case B – Proposed method vs. Mixed-code GA with search enhancement**

| Problem Instance | Proposed method | | Mixed-code GA with search enhancement (Zhang et al. 2015) | |
|---|---|---|---|---|
| | **Optimization Results** | **No. of Lambert routine calls** | **Optimization Results** | **No. of Lambert routine calls** |
| Case B-1 | **q** = (1, 2, 3, 4)<br>*J* = 6.36 *km/s* | 3,526,933 | **q** = (1, 2, 3, 4)<br>*J* = 6.40 *km/s* | 960,000 |
| Case B-2 | **q** = (6, 8, 5, 1, 4, 3, 7, 2)<br>*J* = 16.41 *km/s* | 25,392,677 | **q** = (8, 5, 1, 4, 7, 3, 2, 6)<br>*J* = 19.15 *km/s* | 7,680,000 |
| Case B-3 | **q** = (12, 11, 1, 5, 13, 3, 7, 2, 4, 15, 16, 6, 14, 8, 9, 10)<br>*J* = 37.33 *km/s* | 182,852,203 | **q** = (11, 9, 1, 4, 12, 10, 3, 7, 2, 5, 6, 16, 14, 8, 15, 13)<br>*J* = 57.98 *km/s* | 61,440,000 |

# VI. Conclusions

A novel two-phase framework to solve the optimal multi-target Lambert rendezvous problem is proposed in this paper. Elementary solutions of single target problems associated with all departure-arrival pairs are prepared in the first phase. The second phase search for the best rendezvous sequence and trajectories by solving a variant of TSP formulated using the elementary solutions. Case studies to solve single-target and multi-target rendezvous problems were conducted to demonstrate the effectiveness of the proposed framework.

Additional studies on the enhancing the performance of the proposed framework and systematic comparison with other methodologies can be considered as potential future research. Improvement of Phase A (single-target rendezvous) by considering gravitational perturbation (e.g. $J_2$ term) and its implementation within the proposed framework is promising future work, as well. Consideration of the target selection out of a larger candidate objects with the proposed multi-target rendezvous problem and development of a new framework that can simultaneously handle the combined problem can be another interesting subject for further study.

## References

Ahn, J., and Lee, S., 2013. Lambert Algorithm Using Analytic Gradients. Journal of Guidance, Control, and Dynamics 36, 1751-1761.

Ahn, J., Lee, S., and Bang, J., 2015. Acceleration of zero-revolution Lambert's algorithms using table-based initialization. Journal of Guidance, Control, and Dynamics 38, 335-342.

Alfriend, K. T., Lee, D., and Creamer, N. G., 2005. Optimal Servicing of Geosynchronous Satellites," Journal of Guidance, Control, and Dynamics 29, 203–206.

Bang, J., and Ahn, J., 2016. Optimal Multi-Target Lambert Rendezvous. In: AIAA/AAS Astrodynamics Specialist Conference, (Long Beach, CA), Sept.

Barbee, B. W., Alfano, S., Pinon, E., Gold, K., and Gaylor, D., 2012. Design of spacecraft missions to remove multiple orbital debris objects. Advances in the Astronautical Sciences 144, 93–110.

Bérend, N., and Olive, X., 2016. Bi-objective optimization of a multiple-target active debris removal mission. Acta Astronautica 122, 324–335.

Battin, R. H., 1999. An Introduction to the Mathematics and Methods of Astrodynamics, AIAA, Reston, VA.

Cerf, M., 2013. Multiple Space Debris Collecting Mission-Debris Selection and Trajectory Optimization. Journal of Optimization Theory and Applications 156, 761–796.

Cerf, M., 2015. Multiple Space Debris Collecting Mission: Optimal Mission Planning. Journal of Optimization Theory and Applications 167, 195–218.

Çetinkaya, C., Karaoglan, I., and Gökçen, H., 2013. Two-stage vehicle routing problem with arc time windows: A mixed integer programming formulation and a heuristic approach. European Journal of Operational Research 230, 539–550.

Chen, T., van Kampen, E., Yu, H., and Chu, Q. P., 2013. Optimization of Time-Open Constrained Lambert Rendezvous Using Interval Analysis. Journal of Guidance, Control, and Dynamics 36, 175–184.

Chilan, C. M., and Conway, B. A., 2013. Automated Design of Multiphase Space Missions Using Hybrid Optimal Control. Journal of Guidance, Control, and Dynamics 36, 1410–1424.

Conway, B. A., Chilan, C. M., and Wall, B. J., 2007. Evolutionary principles applied to mission planning problems. Celestial Mechanics and Dynamical Astronomy 97, 73–86.

Englander, J. A., Conway, B. A., and Williams, T., 2012. Automated Mission Planning via Evolutionary Algorithms. Journal of Guidance, Control, and Dynamics 35, 1878–1887.

Garaix, T., Artigues, C., Feillet, D., and Josselin, D., 2010. Vehicle routing problems with alternative paths: An application to on-demand transportation. European Journal of Operational Research 204, 62–75.

Gooding, R. H., 1990. A procedure for the solution of Lambert's orbital boundary-value problem. Celestial Mechanics and Dynamical Astronomy 48, 1990, 145–165.

Izzo, D., Vinkó, T., Bombardelli, C., Brendelberger, S., and Centuori, S., 2007. Automated asteroid selection for a 'grand tour' mission. In: Proceedings of 58th International Astronautical Congress, 4603–4610.

Izzo, D., Getzner, I., Hennes, D., and Simões, L. F., 2015. Evolving Solutions to TSP Variants for Active Space Debris Removal. In: Genetic and Evolutionary Computation Conference (GECCO 2015), (Madrid, Spain), July.

Lee, S., Ahn, J., and Bang, J., 2016. Dynamic selection of zero-revolution Lambert algorithms using performance comparison map. Aerospace Science and Technology 51, 96-105.

Prussing, J. E., 2000. A class of optimal two-impulse rendezvous using multiple-revolution Lambert solutions. Journal of Astronautical Sciences 48, 131–148.

Ross, I. M., and D'Souza, C. N., 2005. Hybrid Optimal Control Framework for Mission Planning. Journal of Guidance Control and Dynamics 28, 686–697.

Shen, H., and Tsiotras, P., Optimal Two-Impulse Rendezvous Using Multiple-Revolution Lambert Solutions. Journal of Guidance Control and Dynamics 26, 50–61.

Ticha, H. B., Absi, N., Feillet, D., and Quilliot, A., 2017. Empirical analysis for the VRPTW with a multigraph representation for the road network. Computers and Operations Research 88, 103–116.

Wall, B. J., and Conway, B. A., 2009. Genetic algorithms applied to the solution of hybrid optimal control problems in astrodynamics. Journal of Global Optimization 44, 493–508.

Yu, J., Chen, X. Q., and Chen, L. H., 2015. Optimal planning of LEO active debris removal based on hybrid optimal control theory. Advances in Space Research 55, 2628–2640.

Zhang, G., Zhou, D., and Mortari, D., 2011. Optimal two-impulse rendezvous using constrained multiple-revolution Lambert solutions. Celestial Mechanics and Dynamical Astronomy 110, 305–317.

Zhang, J., Parks, G. T., Luo, Y., and Tang, G., 2014. Multispacecraft Refueling Optimization Considering the J2 Perturbation and Window Constraints. Journal of Guidance, Control, and Dynamics 37, 111–122.

Zhang, J., Luo, Y., Li, H., and Tang, G., 2015. Analysis of Multiple Asteroids Rendezvous Optimization using Genetic Algorithms. In: Proceedings of 2015 IEEE Congress on Evolutionary Computation (CEC), 596–602.